\documentclass{article}
\usepackage{amssymb}
\usepackage{dsfont}
\usepackage{textcomp}
\usepackage{amsmath}
\usepackage{txfonts}
\usepackage[colorlinks=true]{hyperref}
\usepackage{graphicx}
\usepackage[T1]{fontenc}
\usepackage{stmaryrd}
\usepackage{mathrsfs}

\title{The exponential-logarithmic equivalence classes of surreal numbers.}
\author{Salma Kuhlmann and Micka\"{e}l Matusinski}

\newtheorem{thm}{Theorem}[section]    
\newtheorem{lem}[thm]{Lemme}          
\newtheorem{prop}[thm]{Proposition}    

\newtheorem{conj}[thm]{Conjecture}
\newtheorem{definition}[thm]{Definition}
\newenvironment{defn}{\begin{definition}\rm}{\end{definition}}

\newtheorem{remarque}[thm]{Remark}
\newenvironment{rem}{\begin{remarque}\rm}{\end{remarque}}

\newtheorem{nota}[thm]{Notation}

\newtheorem{prf}{\textit{Proof.}}

\newenvironment{demo}{\begin{prf}\rm}{\hfill$\Box$\end{prf}}

\newtheorem{exemple}[thm]{Example}
\newenvironment{ex}{\begin{exemple}\rm}{\end{exemple}}

\newcommand{\n}{\par\noindent}

\begin{document}
\maketitle

\section{Introduction.}
In his monograph \cite{gonshor_surreal}, H. Gonshor showed that Conway's real closed field of surreal numbers carries an exponential and logarithmic map. Subsequently, L. van den Dries and P. Ehrlich showed in \cite{ehrlich-vdd:surreal-exp} that it is a model of the elementary theory of the field of real numbers with the exponential function. In this paper, we give a complete description of the exponential equivalence classes (see Theorem \ref{theo:kappa-equivEL}) in the spirit of the classical Archimedean and multiplicative equivalence classes (see Theorem \ref{theo:NO-series-gene} and Proposition \ref{prop:N0-Hahn}). This description is made in terms of a recursive formula as well as a sign sequence formula for the family of representatives \emph{of minimal length} of these exponential classes. 

This result can be seen as a step towards proving that the field of surreal numbers \textbf{No} can be described as an exponential-logarithmic series field $\mathds{K}^{\textrm{EL}}$ (for some subfield of generalized series $\mathds{K}$ of \textbf{No}) and a field of transseries $\mathds{T}$. Indeed, we conjecture that our representatives of the exponential classes are the fundamental monomials of $\mathds{K}$ - say the \textbf{initial fundamental monomials} - in the sense of \cite{matu-kuhlm:hardy-deriv-EL-series}, or the \textbf{log-atomic elements} (i.e. the monomials which remains monomials by taking iterated log's) in the sense of \cite{vdh:transs_diff_alg,schm01}. Such a description would allow us to exploit these representations to introduce derivations on the surreals. Indeed, we know how to define derivations on exponential-logarithmic series fields  \cite{matu-kuhlm:hardy-deriv-gener-series, matu-kuhlm:hardy-deriv-EL-series} and on transseries fields \cite{vdh:autom-asymp, schm01}.

In Section \ref{sect:defi}, we give a concise summary of the recursive definitions of the field operations on \textbf{No}, as well as the definition of the exponential $\exp$ and logarithmic $\log$ maps, and generalized epsilon numbers $\epsilon_{\mathrm{\mathbf{No}}}$. Our exposition is based on \cite{gonshor_surreal}. Of particular interest to us is the analysis of certain equivalence relations on the surreal numbers. Conway \cite[p 31-32]{conway_numb-games} introduced and studied the $\omega$-map to give a complete system $\omega^{\mathrm{\mathbf{No}}}$ (:= the image of \textbf{No} under this map) of representatives of the Archimedean additive equivalence relation. In \cite[Theorem 5.3]{gonshor_surreal}, exploiting the convexity of the equivalence classes, Gonshor describes such a representative $\omega^a$ as the \textit{unique} surreal of minimal length in a given class. By a simple modification of their arguments, we describe a complete system $\omega^{\omega^{\mathrm{\mathbf{No}}}}$ of representatives of the Archimedean multiplicative equivalence relation. 
In Section \ref{sect:kappa}, we introduce and study what we call the $\kappa$-map to give a complete system $\kappa_{\mathrm{\mathbf{No}}}$ (:= the image of \textbf{No} under this map) of representatives of the exponential equivalence relation (Definition \ref{defi:equiv-exp} and Theorem \ref{theo:kappa-equivEL}). We observe that:
\begin{center}
$\epsilon_{\textbf{No}}\subsetneq\kappa_{\textbf{No}}\subsetneq \omega^{\omega^{\textbf{No}}} \subsetneq\omega^{\textbf{No}}\subsetneq\textbf{No}$.
\end{center}

Section \ref{sect:kappa-sign} is devoted to establish the sign sequences of these representatives (see Theorem \ref{theo:kappa-signe}). Then, in Section \ref{sect:TEL}, we introduce the notion of Transserial-Exp-Log fields, which unifies the notion of transseries and exp-log series. We conjecture that \textbf{No} is such a TEL field (Conjecture \ref{conj:TEL}).

We thank Joris van der Hoeven : we started with him the project of endowing \textbf{NO} with a derivation, and since then had several valuable exchanges on this topic in connection with the fields of surreal numbers and of 
transseries. This helped to lead us to the results we present here.

\section{Preliminaries on surreal numbers.}\label{sect:defi}
\subsection{Inductive definitions.}
This section is dedicated to fix the notations and recall some of the definitions and results obtained in \cite{conway_numb-games} and \cite{gonshor_surreal}. 
We will also use results from \cite{ehrlich:simpl-hierarc} and \cite{ehrlich-vdd:surreal-exp}.

We denote by  $\textbf{On}$ the proper class of \textbf{ordinal numbers}. A \textbf{surreal number} is any map $a:\alpha\rightarrow \{\ominus,\oplus\}$ where $\alpha\in\textbf{On}$. Any surreal may be canonically represented as a sequence of pluses and minuses called its \textbf{sign sequence}. The ordinal $\alpha$ on which the surreal number $a$ is defined is called its \textbf{length}. Set  $l(a):=\alpha$. The proper class of all surreal numbers is denoted by $\textbf{No}$. 

It is totally ordered as follows. Let $a$ and $b$ be any two surreal numbers  with $l(a)\leq l(b)$. Consider the sign sequence of $a$ completed with 0's so that the two sign sequences have same length. Then consider the \textbf{lexicographical order} between them, denoted by $\leq$, based on the following relation: \begin{center}
$\ominus<0<\oplus$.
\end{center}

A subset of  $\textbf{No}$ is said to be \textbf{initial} if it is made of surreals whose lengths are less than or equal to a given ordinal. As in \cite{ehrlich:simpl-hierarc},  $\textbf{No}$ is endowed with a partial ordering called the \textbf{simplicity }:\begin{center}
 $a\textrm{ is simpler than }b\textrm{, write }a<_sb \Leftrightarrow a\textrm{ is a proper initial segment of }b$.
\end{center} 
Given a surreal  $a\in \textbf{No}$, any surreal number $a'<_sa$ is called a \textbf{truncation} of $a$.

The construction of surreal numbers makes heavy use of the key idea of "cuts between sets" in the classical construction of real numbers. Take two pairs  $(L,R)$ (as "Left" and "Right") and $(L',R')$ of \emph{subsets} of \textbf{No} with $L<R$ and $L'<R'$. $(L',R')$ is said to be \textbf{cofinal} in  $(L,R)$ if for any $(a,b)\in L\times R$, there is $(a',b')\in L'\times R'$, such that $a<a'<b'<b$. Note that cofinality is a transitive relation. The following results are  fundamental. In particular, given a surreal  $a\in\textbf{No}$, we will considere the pair $(L_a,R_a)$ where $L_a=\{b\in\textbf{No}\ ;\ b<a\textrm{ and }b<_sa\}$ and $R_a=\{b\in\textbf{No}\ ;\ b>a\textrm{ and }b<_sa\}$.

\begin{thm}[Existence and cofinality]\label{theo:cofin}
\begin{enumerate}
    \item \cite[Theorem 2.1]{gonshor_surreal} For any pair $(L,R)$ with $L<R$, there exists a \emph{unique} surreal  $a\in \textbf{No}$ of minimal length such that  $L<a<R$. We set $\langle L|R\rangle:=a\in\textbf{No}$, and call it the \textbf{cut} between $L$ and $R$.
\item  \cite[Theorem 2.6]{gonshor_surreal} Suppose that $\langle L\ |\ R\rangle=a\in\textbf{No}$, $L'<a<R'$ for some  pair $(L',R')$ cofinal in  $(L,R)$, then $\langle L'\ |\ R'\rangle=a$. 
\item \cite[Theorems 2.8 and 2.9]{gonshor_surreal} One always has $a=\langle L_a\ |\ R_a\rangle$ and, for any cut $\langle L'\ |\ R'\rangle=a$, then $(L',R')$ is cofinal in $(L_a,R_a)$.
\end{enumerate}
\end{thm}

\begin{defn}
 This representation $a=\langle L_a\ |\ R_a\rangle$ of $a$ is called   \textit{its} \textbf{ canonical cut}.  By abuse of notation, we also denote the canonical cut by $a=\langle a^L\ |\ a^R\rangle$ where $a^L$ and $a^R$ are general elements of the canonical sets  $L_a$ and $R_a$ (e.g. $a^L=n$ if $L_a=\mathbb{N}$).
\end{defn}

The preceding results allow to define maps and operations on surreal numbers \emph{recursively}, that is by induction on the length of the surreals considered since $a^L$ and $a^R$ are simpler than $a$. Concerning the algebraic structure of  $\textbf{No}$, we have:

\begin{thm}\label{theo:NO-reel-clos}
\begin{enumerate}
    \item \cite[Theorems 3.3, 3.6, 3.7 and Ch. 5, Sect. D]{gonshor_surreal} The proper class $\textbf{No}$  endowed with its lexicographical ordering  $\leq$ and the following operations, is a real closed field (in the sense of proper classes): for any $a,b\in\textbf{No}$,
\begin{center}
$\begin{array}
{llcl}\emph{addition: }&a+b&:=&\langle a^L+b,\ a+b^L\ |\ a^R+b,\ a+b^R\rangle\\
 \emph{neutral\ element}: &0&=&\langle\emptyset\ |\ \emptyset\rangle ; \\
\emph{multiplication: }&a.b&:=& \langle a^L.b+a.b^L-a^L.b^L,\ a^R.b+a.b^R-a^R.b^R\ |\ \\ 
&&&\ \ \ \ \ \ a^L.b+a.b^R-a^L.b^R,\ a^R.b+a.b^L-a^R.b^L\rangle  \\
\emph{neutral\ element}: &1&=&\oplus=\langle 0\ |\ \emptyset\rangle. \\
\end{array}$
\end{center}
\item \cite[Theorems 9 and 19]{ehrlich:simpl-hierarc} Any divisible ordered abelian group, respectively any real closed field, is isomorphic to an initial subgroup of  $(\textbf{No},+)$, respectively an initial subfield of  $(\textbf{No},+,.)$.
\end{enumerate}
\end{thm}
Note that Gonshor proves also, for the cited operations as well as for the following maps, the so-called \textbf{uniformity properties} \cite[Theorems 3.2, 3.5 etc.]{gonshor_surreal}. For example in the case of the addition, this means that $a+b$ may be obtained by taking any cuts $a=\langle L\ |\ R\rangle$ and $b=\langle L'\ |\ R'\rangle$ instead of taking their canonical cuts, and applying the same formula as above. In other words, \emph{the formulas do not depend on the cuts for which $a^L,\ a^R,\ b^L,\ b^R$ are general elements}.

Any \textbf{real number} is identified to the sign sequence corresponding to its binary expansion. Thus they are the sequences of length  finite or equal to $\omega$, the later being non ultimately constant. These sequences form a subfield of the surreals having the least upper bound property, i.e. a copy of the ordered set of the reals: $\mathbb{R}\subset \textbf{No}$.

Surreal numbers $a\in\textbf{No}$ such that $a>\mathbb{R}$ or $a<\mathbb{R}$ are called \textbf{infinitely large}: we denote by $\textbf{No}^{\gg 1}$ their proper class. Those that verify $0<a<\mathbb{R}_{>0}$ or $0>a>\mathbb{R}_{<0}$ are called \textbf{infinitesimals}: we denote by $\textbf{No}^{\ll 1}$ their proper class. 

Any ordinal number $\alpha\in\textbf{On}$ is identified to the surreal number whose sign sequence is a sequence of $\oplus$'s of length $\alpha$. The operations precedingly defined on the surreal numbers correspond to the \textbf{natural sum and product} (see e.g. \cite{hausdorff-mengenlehre}). So we consider $\textbf{On}\subset\textbf{No}$.

In the following theorem, we sum up the main results about the so-called \textbf{Conway normal form} of a surreal number:

\begin{thm}[Conway normal form]\label{theo:NO-series-gene}
\begin{enumerate}
    \item \cite[Theorems 5.1 to 5.4]{gonshor_surreal} The recursive formula:
\begin{center}
$\forall a\in\textbf{No}$, $\omega^{a}:=\langle 0,\ n.\omega^{a^L}\ |\ \omega^{a^R}/2^n\}$ 
\end{center}
(where it is understood that $n\in\mathbb{N}$), defines a map:
\begin{center}
$\begin{array}{llcl}
\Omega:&\textbf{No}&\rightarrow&\textbf{No}\\
&a&\mapsto&\Omega(a):=\omega^a
\end{array}$
\end{center} with values in $\textbf{No}_{>0}$ and that extends the exponentiation with base $\omega$ of the ordinals. Moreover, for any $a\in\textbf{No}$, $\omega^a$ is the representative of minimal length of its \textbf{Archimedean equivalence class} ($\forall x,y \in \textbf{No}_{>0},\ x\sim_{arch} y\Leftrightarrow \exists n\in\mathbb{N},\ n.x\geq y \geq x/n$).
\item \cite[Theorems 5.5 to 5.8]{gonshor_surreal} The field of surreals is a field of \textbf{generalized series} \cite[Definition 2.5]{matu-kuhlm:hardy-deriv-gener-series} in the following sense: any surreal number $a\in\textbf{No}$ can be written uniquely as $a=\displaystyle\sum_{i<\lambda}\omega^{a_i}r_i$ where the transfinite sequence $(a_i)_{i<\lambda}$ is strictly decreasing and, for any $i$, $r_i\in\mathbb{R}\setminus\{0\}$. In particular, $\omega^{\textbf{No}}$ is seen as the group of \textbf{(generalized) monomials}: \begin{center}
$\textbf{No}=\mathbb{R}((\omega^{\textbf{No}}))$.
\end{center}
\end{enumerate}
\end{thm}

In particular, the proper class $\omega^{\textbf{No}}$ is a complete system of representatives of the Archimedean equivalence classes of  $\textbf{No}_{>0}$, each of its elements being the reprensentative of minimal length in its class.

Moreover, note that the map denoted by  $Ind$ in \cite{gonshor_surreal}, which sends any surreal $a$ to the exponent $a_0$ of the leading monomial of its normal form, is a \textbf{valuation}. It is the natural valuation of the real closed field $\textbf{No}$ \cite{kuhl:ord-exp}. In particular, $Ind(\omega^{a})=a$ for any $a\in\textbf{No}$.

We deduce immediately that $\textbf{No}$ is a \textbf{Hahn field of series} \cite{hahn:nichtarchim} in the following sense:
\begin{prop}\label{prop:N0-Hahn}
\begin{enumerate}
    \item For any $a\in \textbf{No}$, $\omega^{\omega^a}$ is the representative of minimal length in its \textbf{equivalence class of comparability} ($\forall x,y\in\textbf{No}_{>0}^{\gg 1},\ x\sim_{comp} y\Leftrightarrow (\exists n\in\mathbb{N},\ x^n\geq y \geq x^{1/n}$,  and we set $ 1/x\sim_{comp}x$).
\item Any surreal $a\in\textbf{No}$ can be written uniquely $a=\displaystyle\sum_{i<\lambda}\left(\displaystyle\prod_{j<\lambda_i} \left(\omega^{\omega^{b_{i,j}}}\right)^{s_{i,j}}\right)r_i$ where, for any $i$, we identify $\displaystyle\prod_{j<\lambda_i} \left(\omega^{\omega^{b_{i,j}}}\right)^{s_{i,j}}=\omega^{a_i}$ where $a_i=\displaystyle\sum_{j\in\lambda_i}\omega^{b_{i,j}}$ (so the tranfinite sequences $(a_i)_{i<\lambda}$ and for any $i$, $(b_{i,j})_{j<\lambda_i}$, are strictly decreasing, and for any $i,j$, $s_{i,j},r_i\in\mathbb{R}\setminus\{0\}$). In particular, we call $\omega^{\omega^{\textbf{No}}}$ the chain of \textbf{fundamental monomials} of \textbf{No} \cite[Definition 2.2]{matu-kuhlm:hardy-deriv-gener-series}:
\begin{center}
$\textbf{No}=\mathbb{R}\left(\left(\left(\omega^{\omega^{\textbf{No}}} \right)^{\mathbb{R}}\right)\right)$.\\
\end{center}
\end{enumerate}
\end{prop}

The proper class $\omega^{\omega^{\textbf{No}}}$ is a complete set of representatives of the comparability classes of $\textbf{No}_{>0}^{\gg 1}$, each of its elements being the one of minimal length in its class.\\

Following the ideas of Conway, Gonshor described also a proper class of \textbf{generalized epsilon numbers}  \cite[p.35]{conway_numb-games}:
\begin{thm}\label{theo:eps-nbres}
The recursive formula: \begin{center}
$\forall a\in\textbf{No}$, $\epsilon_{a}:=\langle \Omega^n(1),\ \Omega^n(\epsilon_{a^L}+1)\ |\ \Omega^n(\epsilon{a^R}-1)\}$
\end{center}
(where it is understood that $n\in\mathbb{N}$), defines a map:
\begin{center}
$\begin{array}{llcl}
\epsilon: &\textbf{No}&\rightarrow&\textbf{No}\\
a&\mapsto&\epsilon(a):=\epsilon_a
\end{array}$
\end{center} with values in $\textbf{No}_{>0}^{\gg 1}$ bigger than $\Omega^n(1)$ for any $n$, and which extends the classical epsilon numbers map on ordinals \cite{sierpinski_ordinal-cardinal}. Moreover,  $\epsilon(\textbf{No})$ is the proper class of all the fixed points of the map  $\Omega$ : \begin{center}
$\forall a\in\textbf{No},\ \omega^{\epsilon_a}=\epsilon_a$.\\
\end{center} 
\end{thm}

Gonshor defines inductively in \cite{gonshor_surreal} an \textbf{exponential map} $\exp:(\textbf{No},+)\rightarrow(\textbf{No}_{>0},.)$ which is surjective (and consequently its inverse, the \textbf{logarithmic map} $\log:(\textbf{No}_{>0},.)\rightarrow(\textbf{No},+)$) using the Taylor expansion of the real exponential map. Here we give only the results we will need in this article:

\begin{thm}\label{theo:exp}
The exponential map $\exp:\textbf{No}\rightarrow\textbf{No}_{>0}$ and its inverse $\log:\textbf{No}_{>0}\rightarrow\textbf{No}$, coincide on the real numbers with the usual real exponential and logarithmic maps, and are such that:
\begin{center}
$\begin{array}{rlcl}
1.\ \textrm{for any }a\in\textbf{No},& \exp(\displaystyle\sum_{i<\lambda}\omega^{a_i}r_i)&=& \omega^{\sum_{i<\lambda}\omega^{g(a_i)}r_i}\\
\textrm{where }& g(a)&:=&\langle  Ind(a),\ g(a^L)\ |\ g(a^R)\rangle ;\\
2.\ \textrm{for any }b\in\textbf{No},& \log(\displaystyle\prod_{j<\lambda} \left(\omega^{\omega^{b_{j}}}\right)^{s_{j}})&=& \displaystyle\sum_{j<\lambda}\omega^{h(b_j)}\\
\textrm{where }& h(b)&:=&\langle  0,\ h(b^L)\ |\ h(b^R),\ \omega^b/2^n\rangle.
\end{array}$
\end{center}
We have $\exp=\log^{-1}$ and $h=g^{-1}$ on $\textbf{No}$.
\end{thm}
Note that $g$ is not defined at 0.

\subsection{Sign sequences.}

Gonshor obtains in \cite{gonshor_surreal} detailed results on the sign sequences of the surreals under the various operations and maps. We will use repeatedly his results. Moreover, we introduce in this section some new operations on the sign sequences.

\begin{defn}
\begin{itemize}
    \item Given two surreal numbers $a,b\in\textbf{No}$, we define their \textbf{concatenation} $a\binampersand b$ as the juxtaposition of their sign sequences. We note that $l(a\binampersand b)=l(a)+l(b)$, the \textsl{ordinal sum} of $l(a)$ and $l(b)$.

\item As in \cite[Theorem 9.5]{gonshor_surreal}, for any surreal number  $a\in\textbf{No}$, we can write its sign sequence as the following transfinite concatenation:
\begin{center}
$a=\alpha_0\oplus\binampersand \beta_0\ominus\binampersand \alpha_1\oplus\binampersand \beta_1\ominus\binampersand\cdots$
\end{center}
where for any $\mu\in\textbf{On}$ used in this writing of $a$, we have $\alpha_\mu,\ \beta_\mu\in\textbf{On}$ with in particular $\alpha_\mu$ possibly 0 for $\mu=0$ or for any \textsl{limit} ordinal $\mu$.

\item Given a surreal number $a\in\textbf{No}$, we denote by $a^+$ the total number of pluses in the sign sequence of $a$ ($a^+\in\textbf{No}$). Therefore is the ordinal sum of the packages of   $\oplus$'s:
$$a^+=\displaystyle\sum_{\mu}\alpha_\mu= \displaystyle\binampersand_{\mu}(\alpha_\mu \oplus).$$

\item Given a positive surreal number $a\in\textbf{No}_{>0}$, we derive from it  the surreal number $\flat a$ obtained by suppressing the first $\oplus$ in the sign sequence of $a$.

Given a negative surreal number $a\in\textbf{No}_{<0}$, we derive from it  the surreal number $\sharp a$ obtained by suppressing the first $\ominus$ in the sign sequence of $a$.

\item Given a surreal number in normal form $\displaystyle\sum_{i<\lambda}\omega^{a_i}.r_i$, we define its corresponding \textbf{reduced sequence}  $(a_i^0)_{i<\lambda}$ as follows. For any $i<\lambda$, $a_i^0$ is derived from $a_i$ by suppressing in its sign sequence the following $\ominus$'s:
\begin{itemize}
    \item given an ordinal number $\nu\in\textbf{On}$, if $a_i(\nu)=\ominus$ and if there exists  $j<i$ such that  $a_j(\xi)=a_i(\xi)$ for any  $\xi\leq\nu$, then suppress the $\nu$th $\ominus$;
\item if $i$ is a successor ordinal, $a_{i-1}\binampersand\ominus$ is a truncation of $a_i$ and if $r_{i-1}$ is not a dyadic rationnal number, then suppress the $\ominus$ coming after $a_{i-1}$ in the writing of $a_i$.
\end{itemize}
\end{itemize} 
\end{defn}

Concerning the exponentiation with base $\omega$, we recall \cite[Theorem 5.12 and Corollary 5.1]{gonshor_surreal}:
\begin{thm}\label{theo:signe-omega}
\begin{itemize}
    \item Given a surreal number $a=\alpha_0\oplus\binampersand \beta_0\ominus\binampersand \alpha_1\oplus\binampersand \beta_1\ominus\binampersand\cdots$, for any ordinal $\mu\in\textbf{On}$ intervening in the writing of $a$, we set $\gamma_\mu:=\displaystyle\sum_{\lambda\leq\mu}\alpha_\lambda$ (ordinal sum). Then the sign sequence of $\omega^a$ is:
\begin{center}
$\omega^a=\omega^{\gamma_0}\oplus\binampersand \omega^{\gamma_0+1}\beta_0\ominus\binampersand \omega^{\gamma_1}\oplus\binampersand \omega^{\gamma_1+1}\beta_1\ominus\binampersand\cdots$
\end{center}
\item Given a positive real number $r=\rho_0\oplus\binampersand\sigma_0\ominus\binampersand\cdots$, the sign sequence of $\omega^a.r$ is 
\begin{center}
$\omega^a.r=\omega^a\binampersand \omega^{a^+}\flat\rho_0\oplus\binampersand\omega^{a^+}\sigma_0\ominus\binampersand \omega^{a^+}\rho_1\oplus\binampersand\omega^{a^+}\sigma_1\binampersand\cdots$  
\end{center}(ordinal multiplication).
If $r$ is negative, reverse all the signs in the preceding sequence.

\item Given a surreal number in normal form $\displaystyle\sum_{i<\lambda}\omega^{a_i}.r_i$, its sign sequence is:
\begin{center}
$\displaystyle\binampersand_{i<\lambda}\omega^{a_i^0}.r_i$
\end{center}
where $(a_i^0)_{i<\lambda}$ is the corresponding reduced sequence.\\
\end{itemize}
\end{thm}

\subsection{Generalized epsilon numbers}
We recall \cite[Theorems 9.5 and 9.6]{gonshor_surreal}:

\begin{thm}\label{theo:sign-seq-epsilon}
 \begin{itemize}
    \item A surreal number $a=\alpha_0\oplus\binampersand \beta_0\ominus\binampersand \alpha_1\oplus\binampersand \beta_1\ominus\binampersand\cdots$ is an epsilon number if and only if $\alpha_0\neq 0$,  all $\alpha_\mu $ different from  0  are ordinary epsilon numbers  satisfying $\alpha_\mu  > l.u.b.\{ \alpha_\lambda\ |\ \lambda<\mu\}$ and furthermore $\beta_\mu$  is a multiple of $\omega^{\alpha_\mu \omega}$  for $\alpha_\mu \neq 0$ and a multiple of $\omega^{\gamma_\mu \omega}$    where  $\delta_\mu:=\displaystyle\sum_{\lambda<\mu}\alpha_\lambda$ (ordinal sum) for $\alpha_\mu  = 0$.
\item  Let $\gamma_\mu:=\displaystyle\sum_{\lambda\leq\mu}\alpha_\lambda$ (ordinal sum).  Then
the  $\mu$th block of pluses in  $\epsilon_a$   consists of $\epsilon_{\gamma_\mu}$   pluses and the  $\mu$th
block of minuses of  $(\epsilon_{\gamma_\mu})^\omega \beta_\mu$  minuses.
\end{itemize}
\end{thm}

\section{The kappa map.}\label{sect:kappa}

In this section we define and describe a new map  $\kappa$, which takes naturally place between the Conway-Gonshor maps $\Omega$ and $\epsilon$ (see Remark \ref{rem:gonshor}).
As in \cite[Remark 3.20]{kuhl:ord-exp}, we introduce the notion of \textbf{exponential equivalence relation} for surreal numbers.  Set $\exp^n$ and $\log^n$ for the $n$th iterate of the corresponding maps. 

\begin{defn}\label{defi:equiv-exp}
We set :
\begin{itemize}
    \item the \textbf{exponential equivalence relation} to be: \begin{center}
$\forall x,y\in\textbf{NO}^{\gg 1},\ x\sim_{\exp}y\Leftrightarrow \exists n\in\mathbb{N},\ \log^n(|x|)\leq |y|\leq\exp^n(|x|)$;\\
\end{center}
\item the \textbf{exponential comparison relation} to be:
\begin{center}
 $\forall x,y\in\textbf{NO}^{\gg 1},\ x\gg _{\exp}y\Leftrightarrow \forall n,\ \log^n(|x|)>|y|$.
\end{center}
\end{itemize}
\end{defn}

\begin{thm}\label{theo:kappa-defi}
The recursive formula
$$\forall a\in\textbf{No},\ \kappa(a)=\kappa_a:=\langle \exp^n(0),\ \exp^n (\kappa_{a^L})\ |\ \log^n (\kappa_{a^R})\rangle$$
(where it is understood that $n\in\mathbb{N}$) defines a map
\begin{center}
$\begin{array}{llcl}
\kappa: &\textbf{No}&\rightarrow&\textbf{No}\\
&a&\mapsto&\kappa(a):=\kappa_a
\end{array}$
\end{center}
with values in $\textbf{No}_{>0}^{\gg 1}$  and such that:
\begin{description}
    \item[(i)] for any $a,b\in\textbf{No}$, $a<b\Rightarrow\kappa_a\ll_{\exp}\kappa_b$;
\item [(ii)] there is a uniformity property for this formula (i.e. the recursive formula does not depend on the choice of the cut for $a$).
\end{description}
\end{thm}

\begin{demo} We proceed by transfinite induction on the length of the surreals considered. For $l(a)=0$, i.e. $a=0=\langle \emptyset|\emptyset\rangle$, we have: \begin{center}
$\begin{array}{lcl}
\kappa_0&:=&\langle \exp^n(0)\ |\ \emptyset\rangle\\
&=&\langle n\ |\ \emptyset\rangle\textrm{ by cofinality}\\
&=&\omega.\\
\end{array}$
\end{center}

We consider $a\in\textbf{No}$ with $l(a)>0$, and suppose that the theorem holds for any $b\in\textbf{No}$ with $l(b)<l(a)$. Consider a canonical representation $a=\langle a^L\ |\ a^R\rangle$ of $a$. By the induction hypothesis, since $\kappa_{a^L}>\mathbb{R}$ and $\kappa_{a^L}\ll_{\exp}\kappa_{a^R}$, for any $n\in\mathbb{N}$ we have  $\exp^n(0)<\exp^n (\kappa_{a^L})<\log^n (\kappa_{a^R})$. Thus the recursive formula for  $\kappa_a$ is well defined. Moreover, since $\exp^n(0)$ for $n\in\mathbb{N}$ is cofinal in $\mathbb{R}$ and is always part of the lower elements in the recursive formula, we have $\kappa_a>\mathbb{R}$, meaning that $\kappa_a\in\textbf{No}_{>0}^{\gg 1}$ for any $a\in\textbf{No}$.

(i) To show the property $a<b\Rightarrow\kappa_a\ll_{\exp}\kappa_b$, we proceed as in \cite[Theorems 3.1, 3.4 et 5.2]{gonshor_surreal} by induction on the natural sum of $l(a)$ and $l(b)$. So we consider $a<b$ in $\textbf{No}$. We denote by $c$ the common initial segment of $a$ and $b$. There are two cases. Either $c=a\textrm{ or }b$: in this case, the property follows directly from the recursive definition of $\kappa$, since $a=b^L$ or $b=a^R$ respectively. Or we have $a<c<b$ with $c<_sa$ and $c<_sb$. Then the natural sum of $a$ and $c$, and of $b$ and $c$, is less than the one of $a$ and $b$. By the induction hypothesis, we obtain that $\kappa_a\ll_{\exp}\kappa_c$ and $\kappa_c<_{\exp}\kappa_b$, so $\kappa_a\ll_{\exp}\kappa_b$.

(ii) Concerning the uniformity property, we refer to \cite[Theorems 3.2, 3.5, Corollary to Theorem 5.2 etc]{gonshor_surreal}. Let $a=\langle L\ |\ R\rangle\in\textbf{No}$. By Theorem \ref{theo:cofin} (3.), $(L,R)$ is cofinal in the canonical cut $(L_a,R_a)$ of $a$. By the preceding property, it implies that $(\kappa(L),\kappa(R))$ and $(\kappa(L_a),\kappa(R_a))$ are cofinals, and that $\kappa(L)\ll_{\exp}\kappa_a\ll_{\exp}\kappa(R)\Rightarrow\kappa(L)<\kappa_a<\kappa(R)$ (these inequalities concern elements in $\textbf{No}_{>0}^{\gg 1}$). So we can apply Theorem \ref{theo:cofin} (2.) to obtain 
$\kappa_a=\langle\kappa(L_a)\ |\ \kappa(R_a)\rangle = \langle\kappa(L)\ |\ \kappa(R)\rangle$.
\end{demo}

\begin{ex} We compute some particular values:
\begin{center}
$\begin{array}{lcl}
\kappa_1&:=&\langle \exp^n(0),\ \exp^n(\omega)\ |\ \emptyset\rangle\\
&=&\langle \Omega^n(\omega)\ |\ \emptyset\rangle\textrm{ by cofinality}\\
&=&\epsilon_0\ ;\\
 \kappa_{-1}&:=&\langle \exp^n(0),\ |\ \log^n(\omega) \rangle\\
&=&\langle n\ |\ \omega^{\omega^{-n}}\rangle\textrm{ by cofinality}\\
&=&\omega^{\omega^{-\omega}}\ ;\\
\end{array}$
\end{center} 
\end{ex}

\begin{thm}\label{theo:kappa-equivEL}
A surreal number $b\in\textbf{No}$ is of the form $\kappa_a$ for some $a\in\textbf{No}$ if and only if $b$ is the representative of minimal length in an exponential equivalence class.
\end{thm}
\begin{demo} We got inspired by \cite[Theorem 5.3]{gonshor_surreal}.\\
Consider $a,b\in\textbf{No}$ such that $b\sim_{\exp}\kappa_a\Leftrightarrow  \exists n\in\mathbb{N},\ \log^n(x)\leq y\leq\exp^n(x)$. So we have $\exp^n(0)\ll_{\exp} b$ and $\exp^n (\kappa_{a^L})\ll_{\exp}b\ll_{\exp} \log^n (\kappa_{a^R})$ for any  $n\in\mathbb{N}$. So $\kappa_a$ is an initial part of $b$.\\
Conversely, we show by transfinite induction on $l(b)\in\textbf{On}$ that for any $b\in\textbf{No}_{>0}^{\gg 1}$, there exists $a\in\textbf{No}$ such that $b\sim_{\exp}\kappa_a$. First, by (ii) of Theorem \ref{theo:kappa-defi}, we  note that $a$ is unique whenever it exists.\\
One has $\omega=\kappa_0$. \\
Consider $b=\langle L_b\ |\ R_b\rangle\in\textbf{No}_{>0}^{\gg 1}$ with $l(b)>\omega$ and suppose that for any $c\in\textbf{No}$ with $l(c)<l(b)$ the desired property holds. So, for any $c\in L_b\cup R_b$, we have $c\sim_{\exp} \kappa_d$ for some $d\in\textbf{No}$. Let $L:=\{d\in\textbf{No}\ |\ \exists c\in L_b,\ c\sim_{\exp}\kappa_d\}$ and $R:=\{d\in\textbf{No}\ |\ \exists c\in R_b,\ c\sim_{\exp}\kappa_d\}$. Suppose that $L\cap R\neq\emptyset$. We consider $d\in L\cap R$. So there exist $c_1\in L_b$ and $c_2\in R_b$ such that  $c_1\sim_{\exp}\kappa_d$ and $c_2\sim_{\exp}\kappa_d$. Since $c_1<c<c_2$, it follows that $c\sim_{\exp}\kappa_d$.\\
Suppose now that $L\cap R=\emptyset$. We have $L<R$ (if not we would have $d_1>d_2\Rightarrow\kappa_{d_1}\gg_{\exp}\kappa_{d_2}\Rightarrow c_1>c_2$ with $c_1\in L_b$ and $c_2\in R_b$). Note that, by definition, $\kappa(L)$, respectively $\kappa(R)$, is a complete set of representatives of the exponential-logarithmic equivalence classes containing the elements of $L_b\setminus\{0\}$, respectively of $R_b$. There are 3 different cases:\\
$\bullet$ either there are $a\in L$ and $n\in\mathbb{N}$ such that $\exp^n(\kappa_a)\geq b$. Consider $c_1\in L_b$ such that  $c_1\sim_{\exp}\kappa_a$. Thus $c_1<b\leq\exp^n(\kappa_a)$ and $c_1\sim_{\exp}\kappa_a\sim_{\exp}\exp^n(\kappa_d)$. So $b\sim_{\exp}\kappa_a$;\\
$\bullet$ either there are $a\in R$ and $n\in\mathbb{N}$ such that $\log^n(\kappa_a)\leq b$. Simingly, we obtain $b\sim_{\exp}\kappa_a$ ;\\
$\bullet$ or for any $d_1\in L$, $d_2\in R$ and $n\in\mathbb{N}$, we have $\exp^n(\kappa_{d_1})< b<\log^n(\kappa_{d_2})$. Consider $c_1\in L_b\setminus\{0\}$. There is $d_1\in L$ such that $c_1\sim_{\exp}\kappa_{d_1}$. In particular, $c_1\leq \exp^n(\kappa_{d_1})$ for some $n$. As before, for any $c_2\in R_b$, there exists $d_2\in R$ such that $c_2\geq \log^m(\kappa_{d_2})$. Theorem \ref{theo:cofin} (2.) applies, so $b=\langle \exp^n(0), \exp^n(\kappa_{d_1})\ |\ \log^n(\kappa_{d_2})\rangle=\kappa_a$ where we set $a:=\langle L\ |\ R\rangle$. 
\end{demo}

\section{Sign sequences formulae for the kappa map.}\label{sect:kappa-sign}
For any $a\in\textbf{No}$ and $n\in\mathbb{N}$, we denote by $\kappa_{a,n}:=\log^n(\kappa_a)$ and $\kappa_{a,-n}=\exp^n(\kappa_a)$, where $\log^0$ and $\exp^0$ are equal to the identity map on $\textbf{No}$. The aim of this section is to show that $\kappa_{a,n}$ and $\kappa_{a,-n}$ for any $n$  are elements of $\omega^{\omega^{No}}$, and to give their sign sequence. In order to do so, we introduce the following auxiliary map $\iota$ and give its sign sequence. 

\begin{nota}\label{nota:gamma}
\begin{itemize}
\item Consider a surreal number $a=\alpha_0\oplus\binampersand \beta_0\ominus\binampersand \alpha_1\oplus\binampersand \beta_1\ominus\binampersand\cdots$. For any  ordinal $\mu\in\textbf{On}$ intervening in the writing of $a$, we set: \begin{center}
$\gamma_\mu:=\displaystyle\sum_{\lambda\leq\mu}\alpha_\lambda\ \  \textrm{ (ordinal sum).}$
\end{center}
\item We set $\Omega^0$ to be the constant map equal to 1 in $\textbf{No}$ and  $\epsilon_{\flat 0}:=0$.
\end{itemize}
\end{nota}

\begin{lem}\label{lemme:iota}
 The recursive formula:
\begin{center}
$\forall a\in\textbf{No},\ \iota(a)=\iota_a:=\langle \iota_{a^L}\binampersand\Omega^n(\epsilon_{\flat(a^L)^+}+1)\oplus\ |\ \iota_{a^R}\binampersand n\ominus\rangle$.
\end{center}
(where it is understood that $n\in\mathbb{N}$) defines a map
\begin{center}
$\begin{array}{llcl}
\iota: &\textbf{No}&\rightarrow&\textbf{No}\\
a&\mapsto&\iota(a):=\iota_a
\end{array}$
\end{center}
such that for any $a=\alpha_0\oplus\binampersand \beta_0\ominus\binampersand \alpha_1\oplus\binampersand \beta_1\ominus\binampersand\cdots$, the sign sequence of $\iota_a$ is given by:
\begin{center}
$\iota_a=\epsilon_{\flat\gamma_{0}}\oplus\binampersand\omega\beta_0\ominus \binampersand\epsilon_{\flat\gamma_1}\oplus\binampersand\omega\beta_1\ominus \binampersand\cdots$.
\end{center}
A uniformity property holds for this map.
\end{lem}
\begin{demo} 
We proceed by transfinite induction on  $l(a)\in\textbf{On}$. We have $\iota_0=\langle\emptyset\ |\ \emptyset\rangle=0$.\\
Given a surreal number  $a\in\textbf{No}$ with $l(a)>0$, we suppose that for any $b$ with $l(b)<l(a)$, $\iota_b$ is well defined, with the corresponding sign sequence. Consider a canonical representation $a=\langle a^L\ |\ a^R\rangle$ of $a$. For instance, in the case where $l(a)$ is a successor ordinal and $a=a^L\binampersand\oplus$. So $a^R<_sa^L$ et $a^R>a^L$. More precisely, either $a^L$ is only made of $\oplus$'s in which case $a^R=\emptyset$ (case (i)), or $a^R$ is a truncation of $a^L$, ending with a group of $\ominus$'s of $a^L$ of smaller number than the corresponding one in $a^L$ (case (ii)). In the case (i), the sign sequence formula of $\iota_a$ is:
\begin{center}
$\begin{array}{lcl}
\iota_a&=&\langle \iota_{a^L}\binampersand\Omega^n(\epsilon_{\flat(a^L)^+}+1)\ |\ \emptyset\rangle\\
&=&\langle \epsilon_{\flat(a^L)^+}\binampersand\Omega^n(\epsilon_{\flat(a^L)^+}+1)\ |\ \emptyset\rangle\\
&=&\epsilon_{\flat(a^L)^++1}\\
&=&\epsilon_{\flat a^+}
\end{array}$
\end{center}
In the case (ii), denote by $a^L=\alpha_0\oplus\binampersand \beta_0\ominus\binampersand\cdots\binampersand \alpha_\mu\oplus\binampersand \beta_\mu\ominus\binampersand\cdots$. So $a^R=\alpha_0\oplus\binampersand \beta_0\ominus\binampersand\cdots\binampersand \alpha_\mu\oplus\binampersand \tilde{\beta}_\mu\ominus$ with $\tilde{\beta}_\mu<\beta_\mu$ (in $\textbf{On}$) for some $\mu\in\textbf{On}$. By the induction hypothesis, $\iota_{a^L}=\epsilon_{\flat\gamma_{0}}\oplus\binampersand\omega\beta_0\ominus \binampersand\cdots\binampersand\epsilon_{\flat\gamma_\mu}\oplus \binampersand\omega\beta_\mu\ominus\binampersand\cdots<\epsilon_{\flat\gamma_{0}}\oplus\binampersand\omega\beta_0\ominus \binampersand\cdots\binampersand\epsilon_{\flat\gamma_\mu}\oplus \binampersand\omega\tilde{\beta}_\mu\ominus=\iota_{a^R}$. So we obtain as desired $\iota_{a^L}\binampersand\Omega^n(\epsilon_{\flat(a^L)^+}+1)<\iota_{a^R}\binampersand n\ominus$ for any $n\in\mathbb{N}$. Moreover,
\begin{center}
$\begin{array}{lcl}
\iota_a&=&\langle \iota_{a^L}\binampersand\Omega^n(\epsilon_{\flat(a^L)^+}+1)\ |\ \iota_{a^R}\binampersand n\ominus\rangle\\
&=&\langle (\epsilon_{\flat\gamma_{0}}\oplus\binampersand\omega\beta_0\ominus \binampersand\cdots\binampersand\epsilon_{\flat\gamma_\mu}\oplus \binampersand\omega\beta_\mu\ominus\binampersand\cdots) \binampersand\Omega^n(\epsilon_{\flat(a^L)^+}+1)\ |\ \\
&&\ \ \ \ \ (\epsilon_{\flat\gamma_{0}}\oplus\binampersand\omega\beta_0\ominus \binampersand\cdots\binampersand\epsilon_{\flat\gamma_\mu}\oplus \binampersand\omega\tilde{\beta}_\mu\ominus)\binampersand n\ominus\rangle\\
&=&\epsilon_{\flat\gamma_{0}}\oplus\binampersand\omega\beta_0\ominus \binampersand\cdots\binampersand\epsilon_{\flat\gamma_\mu}\oplus \binampersand\omega\beta_\mu\ominus\binampersand\cdots) \binampersand\epsilon_{\flat(a^L)^++1}\\
&=&\epsilon_{\flat\gamma_{0}}\oplus\binampersand\omega\beta_0\ominus \binampersand\cdots\binampersand\epsilon_{\flat\gamma_\mu}\oplus \binampersand\omega\beta_\mu\ominus\binampersand\cdots) \binampersand\epsilon_{\flat(a)^+}.
\end{array}$
\end{center}
The other cases for $l(a)$ and $a$ are analogous. They are left to the reader as an exercise. \\
Concerning the uniformity property, by the sign sequence formula for $\iota$, we note that for any surreal numbers $a,b\in\textbf{No}$, $a<b\Rightarrow\iota_a<\iota_b$. So the uniformity property  follows from Theorem \ref{theo:cofin} in the same way as in \cite[Theorem 3.2, 3.5, Corollary au Theorem 5.2]{gonshor_surreal} and in Theorem \ref{theo:kappa-defi}.
\end{demo}

The main result of this section is:

\begin{thm}\label{theo:kappa-signe}
\begin{enumerate}
    \item For any $a\in\textbf{No}$, $n\in\mathbb{N}$, we have $\kappa_{a,n},\ \kappa_{a,-n}\in\omega^{\omega^{\textbf{No}}}$. More precisely:
\begin{center}
$\begin{array}
{lcl}
\kappa_{a,n}&=&\omega^{\omega^{b}}\ \textrm{ with }\ b=\iota_a\binampersand n\ominus;\\
\kappa_{a,-n-1}&=&\omega^{\omega^{b}}\ \textrm{ with }\ b=\iota_a\binampersand \Omega^n(\epsilon_{\flat a^+}+1)\oplus.
\end{array}$
\end{center}
In particular, $\kappa_a=\omega^{\omega^{\iota_a}}$.
\item With the Notation \ref{nota:gamma}, for any $a\in\textbf{No}$, $n\in\mathbb{N}$, we have:
\begin{center}
$\begin{array}
{lcl}
\kappa_a&=&\omega^{\omega^{\epsilon_{\flat\gamma_{0}}}}\oplus\binampersand [(\omega^{\omega^{\epsilon_{\flat\gamma_{0}}}}\omega)^2\beta_0]\ominus \binampersand\epsilon_{\flat\gamma_1}\oplus\binampersand [(\epsilon_{\flat\gamma_{1}}\omega)^2\beta_1]\ominus\binampersand\cdots\ ;\\
\kappa_{a,n}&=&\kappa_a\binampersand (\epsilon_{\flat a^+}\omega n)\ominus\ ;\\
\kappa_{a,-n-1}&=&\kappa_a\binampersand (\Omega^{n+2}(\epsilon_{\flat a^+}+1)\oplus.
\end{array}$
\end{center}
\end{enumerate}
\end{thm}

The second point of the theorem follows directly from the first one and the sign sequences for the maps $\Omega$ (Theorem \ref{theo:signe-omega}) and $\iota$ (Lemma \ref{lemme:iota}). The first point will be deduced from the following lemma.

\begin{lem}\label{lemme:h-iota}
For any $a\in\textbf{No}$ and any $\lambda\in\textbf{On}$, if we set $\lambda =\epsilon_\nu +\mu$ where $\epsilon_\nu$ is the unique epsilon number (possibly 0) such that $\epsilon_\nu<\lambda\leq \epsilon_{\nu+1}$, we have:
\begin{center}
$\begin{array}
{lcl}
h(\iota_a\binampersand \lambda\ominus)&=&\omega^{\iota_a\binampersand(\lambda+1)\ominus} ;\\
h(\iota_a\binampersand \lambda\oplus)&=&\omega^{\iota_a}\binampersand \epsilon\oplus\binampersand\flat\mu\oplus\textrm{ avec }\lambda>0.
\end{array}$
\end{center}
\end{lem}
\begin{demo}
We proceed by transfinite induction on $(l(a),\lambda)\in \overrightarrow{\textbf{On}\times\textbf{On}}$ (the lexicographical product of $\textbf{On}$ with itself).\\
Suppose that $l(a)=0$ i.e. $a=0$. By the definition of the map $h$ (see Theorem \ref{theo:exp} (2.)), we have $h(0)=\langle 0\ |\ 1/2^n\rangle=\omega^{-1}$.\\

Moreover, $1=\oplus=\langle 0\ |\ \emptyset\rangle$. So:
\begin{center}
$\begin{array}
{lcl}
h(\oplus)&=&\langle 0,h(0)\ |\ \omega/2^n\rangle\\
&=& \langle \omega^{-1}\ |\ \omega/2^n\rangle\\
&=&\langle \oplus\binampersand\omega \ominus\ |\ \omega\oplus\binampersand\n\ominus\rangle\textrm{ by the Theorem \ref{theo:signe-omega}}\\
&=& \oplus = \omega^0.
\end{array}$
\end{center}

Consider now $(a,\lambda)\in\textbf{No}\times\textbf{On}$ with $\lambda>0$ in $\textbf{On}$. Suppose that the lemma holds for any $(b,\mu)\in\textbf{No}\times\textbf{On}$ with $(l(b),\mu)<_{lex}(l(a),\lambda)$ in $\overrightarrow{\textbf{On}\times\textbf{On}}$. We have $\iota_a\binampersand\lambda\ominus=\langle (\iota_a)^L\ |\ \iota_a\binampersand\mu\ominus\rangle$ with $\mu<\lambda$ in $\textbf{On}$. So, by the definition of $h$ and by the induction hypothesis:
\begin{center}
$\begin{array}{lcl}
h(\iota_a\binampersand\lambda\ominus)&=&\langle 0,\ h((\iota_a)^L)\ |\ h(\iota_a\binampersand\mu\ominus),\ \omega^{\iota_a\binampersand\lambda\ominus}/2^n\rangle\\
&=& \langle 0,\ h(\iota_{a^L}\binampersand\Omega^n(\epsilon_{\flat(a^L)^+}+1)\oplus)\ |\ \omega^{\iota_a\binampersand(\mu+1)\ominus},\ \omega^{\iota_a\binampersand\lambda\ominus}/2^n\rangle\textrm{ (Lemma \ref{lemme:iota})}\\
&=& \langle 0,\ \omega^{\iota_{a^L}}\binampersand\Omega^n(\epsilon_{\flat(a^L)^+}+1)\oplus\ |\ \omega^{\iota_a\binampersand\lambda\ominus}\binampersand(\epsilon_{\flat a+}n)\ominus\rangle\textrm{ by cofinality }\\
&=& \omega^{\iota_a\binampersand\lambda\ominus} \binampersand(\epsilon_{\flat a+}\omega)\ominus\\
&=& \omega^{\iota_a\binampersand(\lambda+1)\ominus}\textrm{ by Theorem \ref{theo:signe-omega}}.
\end{array}$
\end{center}

We set $\lambda=\epsilon_\nu+\mu$ where $\epsilon_\nu<\lambda\leq \epsilon_{\nu+1}$. We have $\iota_a\binampersand\lambda\oplus=\iota_a\binampersand\epsilon_\nu\oplus\binampersand\mu\oplus= \langle\iota_a\binampersand\epsilon_\nu\oplus\binampersand\rho\oplus\ |\  (\iota_a)^R\rangle$ with $\rho<\mu$ in $\textbf{On}$. So, by the definition of $h$ and by the induction hypothesis:
\begin{center}
$\begin{array}{lcl}
h(\iota_a\binampersand\lambda\oplus)&=&\langle 0,\ h(\iota_a\binampersand\epsilon_\nu\oplus\binampersand\rho\oplus)\ |\ h((\iota_a)^R),\ \omega^{\iota_a\binampersand\lambda\oplus}/2^n\rangle\\
&=& \langle 0,\ \omega^{\iota_a}\binampersand \epsilon_\nu\oplus\binampersand\flat\rho\oplus\ |\ h(\iota_{a^R}\binampersand n\ominus),\ \omega^{\iota_a\binampersand\lambda\oplus}/2^n\rangle\\
&=& \langle 0,\ \omega^{\iota_a}\binampersand \epsilon_\nu\oplus\binampersand\flat\rho\oplus\ |\ \omega^{\iota_{a^R}\binampersand(n+1)\ominus},\ \omega^{\iota_a\binampersand\lambda\oplus}/2^n\rangle\\
&&\ \ \textrm{ (note that }\iota_{a^R}\binampersand(n+1)\ominus>\iota_a\binampersand\lambda\oplus)\\
&=& \langle 0,\ \omega^{\iota_a}\binampersand \epsilon_\nu\oplus\binampersand\flat\rho\oplus\ |\ \omega^{\iota_{a^R}\binampersand(n+1)\ominus},\ \omega^{\iota_a}\binampersand(\epsilon_{\flat a^+}\omega^{\lambda})\oplus\binampersand (\epsilon_{\flat a^+}\omega^{\lambda}n)\ominus\rangle\\
&&\ \ \textrm{ (Theorem \ref{theo:signe-omega} and cofinality) }\\
&=& \omega^{\iota_a}\binampersand\epsilon\oplus\binampersand \flat\mu\oplus.
\end{array}$
\end{center}

Consider now $a\in\textbf{No}$ with $l(a)>0$, and suppose that the lemma holds for any $(b,\mu)\in\textbf{No}\times\textbf{On}$ with $l(b)<l(a)$ in $\textbf{On}$. We have $\iota_a=\langle \iota_{a^L}\binampersand\Omega^n(\epsilon_{\flat(a^L)^+}+1)\oplus\ |\ \iota_{a^R}\binampersand n\ominus\rangle$ by Lemma \ref{lemme:iota}. So, by the definition of $h$ and the induction hypothesis, we have:
\begin{center}
$\begin{array}{lcl}
h(\iota_a)&=&\langle 0,\ h((\iota_a)^L)\ |\ h((\iota_a)^R),\ \omega^{\iota_a}/2^n\rangle\\
&=& \langle 0,\ h(\iota_{a^L}\binampersand\Omega^n(\epsilon_{\flat(a^L)^+}+1)\oplus)\ |\ h(\iota_{a^R}\binampersand n\ominus),\ \omega^{\iota_a}/2^n\rangle\textrm{ (Lemma \ref{lemme:iota})}\\
&=& \langle 0,\ \omega^{\iota_{a^L}}\binampersand\Omega^n(\epsilon_{\flat(a^L)^+}+1)\oplus\ |\ \omega^{\iota_{a^R}\binampersand (n+1)\ominus},\ \omega^{\iota_a}\binampersand(\epsilon_{\flat a+}n)\ominus\rangle\textrm{ (Theorem \ref{theo:signe-omega})}\\
&=& \langle 0,\ \omega^{\iota_{a^L}}\binampersand\Omega^n(\epsilon_{\flat(a^L)^+}+1)\oplus\ |\ \omega^{\iota_a}\binampersand(\epsilon_{\flat a+}n)\ominus\rangle\textrm{ by cofinality }\\
&=& \omega^{\iota_a} \binampersand(\epsilon_{\flat a+}\omega)\ominus\\
&=& \omega^{\iota_a\binampersand\ominus}\textrm{ (Theorem \ref{theo:signe-omega})}.
\end{array}$
\end{center}
\end{demo}

\noindent\textbf{\textit{Proof of Theorem \ref{theo:kappa-signe} (1.).}} $\ \ $ We proceed by tranfinite induction on $(l(a),n)\in\overrightarrow{\textbf{On}\times\mathbb{N}}$.

Suppose that $l(a)=0$, i.e. $a =0$. In the proof of the  Theorem \ref{theo:kappa-defi}, we proved that $\kappa_0=\omega$. So $\kappa_0=\omega^{\omega^0}$ as desired.

Let $n\in\mathbb{N}^*$. We suppose that the property holds for any $(b,m)\in\textbf{No}\times\mathbb{N}$ with $(l(b),m)<(l(a),n)$. We have:
\begin{center}
$\begin{array}{lcl}
\kappa_{a,n}&=&\log^n(\kappa_a)\\
&=&\log(\kappa_{a,n-1})\\
&=&\log(\omega^{\omega^{b}})\ \textrm{ with }\ b=\iota_{a}\binampersand(n-1)\ominus\\
&=& \omega^{h(\iota_{a}\binampersand(n-1)\ominus)}\textrm{ (Theorem \ref{theo:exp})}\\
&=& \omega^{\omega^{b}}\ \textrm{ with }\ b=\iota_{a}\binampersand n\ominus\textrm{ (Lemma \ref{lemme:h-iota})}\\
\end{array}$
\end{center}

We also have:
\begin{center}
$\begin{array}{lcl}
\kappa_{a,-1}&=&\exp(\kappa_a)\\
&=&\exp(\omega^{\omega^{\iota_{a}}})\\
&=& \omega^{\omega^{b}}\ \textrm{ with }\ b=g(\omega^{\iota_{a}})\textrm{ (Theorem \ref{theo:exp})}\\
&=& \omega^{\omega^{b}}\ \textrm{ with }\ b=\iota_{a}\binampersand\oplus\\
&&\ \ \textrm{ (Lemma \ref{lemme:h-iota} : }g(\omega^{\iota_{a}})=\iota_{a}\binampersand\oplus\Leftrightarrow \omega^{\iota_{a}}=h(\iota_{a})
\end{array}$
\end{center}
and 
\begin{center}
$\begin{array}{lcl}
\kappa_{a,-n-1}&=&\exp^{n+1}(\kappa_a)\\
&=&\exp(\kappa_{a,-n})\\
&=&\exp(\omega^{\omega^{b}})\ \textrm{ with }\ b=\iota_{a}\binampersand\Omega^{n-1}(\epsilon_{\flat a^+}+1)\oplus\\
&=& \omega^{\omega^{b}}\ \textrm{ with }\ b=g(\omega^{\iota_{a}\binampersand\Omega^{n-1}(\epsilon_{\flat a^+}+1)\oplus}) \textrm{ (Theorem \ref{theo:exp})}\\
&=& \omega^{\omega^{b}}\ \textrm{ with }\ b=\iota_{a}\binampersand\Omega^{n}(\epsilon_{\flat a^+}+1)\oplus \textrm{ (Lemma \ref{lemme:h-iota})}\\
\end{array}$
\end{center}

Consider $a\in\textbf{No}$ with $l(a)>0$, and suppose that the property holds for any $(b,m)\in\textbf{No}\times\mathbb{N}$ with $l(b)<l(a)$ and any $m\in\mathbb{N}$. By the definition of $\kappa$ (Theorem \ref{theo:kappa-defi}) and by the induction hypothesis, we have:
\begin{center}
$\begin{array}{lcl}
\kappa_a&=&\langle \exp^n(0),\ \exp^n (\kappa_{a^L})\ |\ \log^n (\kappa_{a^R})\rangle\\
&=& \langle n,\ \omega^{\omega^{b}}\ |\ \omega^{\omega^{c}},\ \omega^{\iota_a}/2^n\rangle\ \textrm{ with }\ b=\iota_{a^L}\binampersand \Omega^n(\epsilon_{\flat(a^L)^+}+1)\oplus)\textrm{ and }c=\iota_{a^R}\binampersand n\ominus \\
&=&\omega^{\omega^{\iota_{a}}} \textrm{ (by cofinality and definition of }\iota\textrm{, Lemma \ref{lemme:iota})}\\
\end{array}$
\end{center}
This finishes the proof of Theorem \ref{theo:kappa-signe}.\hfill$\Box$

\begin{rem}\label{rem:gonshor}
\begin{enumerate}
\item One can find in \cite{gonshor_surreal} partial results about these sign sequences. Indeed, for any $n\in\mathbb{N}$, any ordinal $\lambda\in\textbf{On}$, any epsilon-number $\epsilon_\mu$ with $\mu\in\textbf{On}$, we have: 
\begin{center}
$\begin{array}{lcl}
\log^n(\omega^{\omega^{\lambda\ominus}})&=&\omega^{\omega^{(\lambda+n)\ominus}}\textrm{ by \cite[Theorem 10.15]{gonshor_surreal}}\\
\exp^n(\omega)&=&\Omega^{n+1}(1)\textrm{ by \cite[Theorem 10.14]{gonshor_surreal}}\\
\exp^n(\epsilon_\mu)&=&\Omega^{n+1}(\epsilon_\mu+1)\textrm{ by \cite[Theorem 10.14]{gonshor_surreal}}\\
\end{array}$.
\end{center}
\item By the Theorems \ref{theo:signe-omega}, \ref{theo:sign-seq-epsilon} and \ref{theo:kappa-signe}, we have:
\begin{center}
$\epsilon_{\textbf{No}}\subsetneq \kappa_{\textbf{No}}\subsetneq\omega^{\omega^{\textbf{No}}} \subsetneq\omega^{\textbf{No}}\subsetneq\textbf{No}$.
\end{center}
\end{enumerate} 
\end{rem}

\section{Transseries and exp-log series.}\label{sect:TEL}

Fields of Transseries \cite{ecalle:dulac, vdh:autom-asymp, vdh:transs_diff_alg} or log-exp series \cite{vdd:LE-pow-series}, and fields of exp-log series \cite{kuhl:ord-exp, kuhl-saharon:kappa-bounded} are important non standard models of the theory of $\mathbb{R}_{\exp}$, which is known to be o-minimal and model complete \cite{wil:modelcomp}. Moreover, these fields can be endowed with derivations that mimic the derivation of germs of real functions in a Hardy field \cite{schm01, dmm:LE-series, matu-kuhlm:hardy-deriv-EL-series}. We propose the following unifying notion for transseries and Exp-Log series, which we believe applies to the field of surreal numbers:

\begin{defn}\label{defi:trans-el}
Let $\mathds{K}=\mathbb{R}((\Gamma))$ be a field of generalized series endowed with a partial logarithm $\log : \Gamma\rightarrow \mathbb{R}((\Gamma^{\succ 1}))$. A complete subfield $\mathbb{L}\subset \mathds{K}$ which contains $\Gamma$ is called a field of \textbf{exp-log-transseries} if the following properties hold:
\begin{description}
    \item[TEL1.] $\textrm{domain }\log=\mathbb{L}_{>0}$.
\item[TEL2.] $\log (\Gamma)= \mathbb{L}^{\succ 1}$.
\item[TEL3.] $\log(1+\epsilon)=\displaystyle\sum_{n=1}^{\infty}\epsilon^n/n \in\mathbb{L}^{\prec 1}$ for any $\epsilon\in \mathbb{L}^{\prec 1}$.
\item[TEL4.] For any sequence of monomials $(m_n)_{n}\subset\Gamma$ such that  for any $n\in\mathbb{N}$, $m_{n+1}\in \textrm{Supp }\log(m_n)$, then there exists a rank $N\in\mathbb{N}$ such that $\log(m_{N+k})=m_{N+k+1}$ for any $k\in\mathbb{N}$.
\end{description}
\end{defn}

Any field of exp-log-transseries is both a field of transseries (our four axioms are specializations of the four axioms of \cite[Definition 2.2.1]{schm01}) and a field of exp-log series. Indeed, denote by $\Phi\subset \Gamma$  the biggest subset of monomials stable by $\log$ (i.e. the set of log-atomic elements in \cite[Section 2.2]{schm01}). Then resuming the notations and terminology of \cite{matu-kuhlm:hardy-deriv-EL-series}, $\Phi$ is a totally ordered set of fundamental monomials. For the corresponding Hahn series field, we set $\mathbb{L}_0=\mathbb{R}((H(\Phi)))\subset \mathbb{L}$ where $H(\Phi)$ denotes the Hahn group over $\Phi$. Then by (TL4) one has that $\mathbb{L}=\mathbb{L}_0^{EL}$.

We believe that:
\begin{conj}\label{conj:TEL}
The field of surreal numbers is an exp-log-transseries field. The surreals $\kappa_{a,n}$ for $a\in\mathbb{No}$, $n\in\mathbb{Z}\}$, are the log-atomic elements as well as the initial fundamental monomials.
\end{conj}

In particular, this would allow us to use the results of \cite{schm01, matu-kuhlm:hardy-deriv-EL-series} to endow the surreal numbers with Hardy type derivations.


\begin{thebibliography}{vdDMM01}

\bibitem[Con01]{conway_numb-games}
J.~H. Conway, \emph{On numbers and games}, second ed., A. K. Peters Ltd.,
  Natick, MA, 2001.

\bibitem[\'E92]{ecalle:dulac}
J.~\'Ecalle, \emph{Introduction aux fonctions analysables et preuve
  constructive de la conjecture de {D}ulac}, Actualit\'es Math\'ematiques.
  [Current Mathematical Topics], Hermann, Paris, 1992.

\bibitem[Ehr01]{ehrlich:simpl-hierarc}
P.~Ehrlich, \emph{Number systems with simplicity hierarchies: a generalization
  of {C}onway's theory of surreal numbers}, J. Symbolic Logic \textbf{66}
  (2001), no.~3, 1231--1258.

\bibitem[Gon86]{gonshor_surreal}
H.~Gonshor, \emph{An introduction to the theory of surreal numbers}, London
  Mathematical Society Lecture Note Series, vol. 110, Cambridge University
  Press, Cambridge, 1986.

\bibitem[Hah07]{hahn:nichtarchim}
H.~Hahn, \emph{{\"U}ber die nichtarchimedischen {G}r\"ossensystem},
  Sitzungsberichte der Kaiserlichen Akademie der Wissenschaften, Mathematisch -
  Naturwissenschaftliche Klasse (Wien) \textbf{116} (1907), no.~Abteilung IIa,
  601--655.

\bibitem[Hau44]{hausdorff-mengenlehre}
F.~Hausdorff, \emph{Mengenlehre}, Dover Publications, New York, N. Y., 1944.

\bibitem[KM11]{matu-kuhlm:hardy-deriv-EL-series}
S.~Kuhlmann and M.~Matusinski, \emph{Hardy type derivations on fields of
  exponential logarithmic series.}, J. Algebra \textbf{345} (2011), 171--189.

\bibitem[KM12]{matu-kuhlm:hardy-deriv-gener-series}
\bysame, \emph{Hardy type derivations in generalized series fields.}, J.
  Algebra \textbf{351} (2012), 185--203.

\bibitem[KS05]{kuhl-saharon:kappa-bounded}
S.~Kuhlmann and S.~Shelah, \emph{{$\kappa$}-bounded exponential-logarithmic
  power series fields}, Ann. Pure Appl. Logic \textbf{136} (2005), no.~3,
  284--296.

\bibitem[Kuh00]{kuhl:ord-exp}
S.~Kuhlmann, \emph{Ordered exponential fields}, Fields Institute Monographs,
  vol.~12, American Mathematical Society, Providence, RI, 2000.

\bibitem[Sch01]{schm01}
M.C. Schmeling, \emph{Corps de transs\'eries}, Ph.D. thesis, Universit\'e
  Paris-VII, 2001.

\bibitem[Sie65]{sierpinski_ordinal-cardinal}
W.~Sierpi{\'n}ski, \emph{Cardinal and ordinal numbers}, Second revised edition.
  Monografie Matematyczne, Vol. 34, Pa\'nstowe Wydawnictwo Naukowe, Warsaw,
  1965.

\bibitem[vdDE01]{ehrlich-vdd:surreal-exp}
L.~van~den Dries and P.~Ehrlich, \emph{Fields of surreal numbers and
  exponentiation}, Fund. Math. \textbf{167} (2001), no.~2, 173--188.

\bibitem[vdDMM97]{vdd:LE-pow-series}
L.~van~den Dries, A.~Macintyre, and D.~Marker, \emph{Logarithmic-exponential
  power series}, J. London Math. Soc. (2) \textbf{56} (1997), no.~3, 417--434.

\bibitem[vdDMM01]{dmm:LE-series}
\bysame, \emph{Logarithmic-exponential series}, Proceedings of the
  {I}nternational {C}onference ``{A}nalyse \& {L}ogique'' ({M}ons, 1997), vol.
  111, 2001, pp.~61--113.

\bibitem[vdH97]{vdh:autom-asymp}
J.~van~der Hoeven, \emph{Asymptotique automatique}, Universit\'e Paris VII,
  Paris, 1997, Th\`ese, Universit\'e Paris VII, Paris, 1997, With an
  introduction and a conclusion in French.

\bibitem[vdH06]{vdh:transs_diff_alg}
\bysame, \emph{Transseries and real differential algebra}, Lecture Notes in
  Mathematics, vol. 1888, Springer-Verlag, Berlin, 2006.

\bibitem[Wil96]{wil:modelcomp}
A.~J. Wilkie, \emph{Model completeness results for expansions of the ordered
  field of real numbers by restricted {P}faffian functions and the exponential
  function}, J. Amer. Math. Soc. \textbf{9} (1996), no.~4, 1051--1094.

\end{thebibliography}
\def\cprime{$'$}
\providecommand{\bysame}{\leavevmode\hbox to3em{\hrulefill}\thinspace}
\providecommand{\MR}{\relax\ifhmode\unskip\space\fi MR }
\providecommand{\MRhref}[2]{%
  \href{http://www.ams.org/mathscinet-getitem?mr=#1}{#2}
}
\providecommand{\href}[2]{#2}

\end{document}